\documentclass[preprint,12pt]{elsarticle}
\usepackage[a4paper, total={6in, 9in}]{geometry}
\usepackage[utf8]{inputenc}
\usepackage{amssymb}
\usepackage{amsmath}
\usepackage{amsfonts}
\usepackage{multirow}
\usepackage{tabularx}
\usepackage{makecell}
\usepackage{comment}
\usepackage{float}
\usepackage{graphicx}
\usepackage{subcaption}
\usepackage{listings}
\usepackage{algorithm}
\usepackage{algpseudocode}
\usepackage{pgfplots}
\usepackage{tabularray}
\usepackage{url}
\UseTblrLibrary{booktabs}
\AfterEndEnvironment{lstlisting}{\leavevmode}

\definecolor{mygreen}{RGB}{0, 153, 0}

\usepackage{todonotes}
\presetkeys{todonotes}{inline,caption={}}{}

\newcommand{\chem}[1]{\ensuremath{\mathrm{#1}}}

\usepackage{amsthm}
\theoremstyle{plain}
\newtheorem{deff}{Definition}

\newtheorem{teo}{Theorem}

\newcommand{\numberset}{\mathbb}

\newcommand{\Z}{\numberset{Z}}

\newcommand{\R}{\numberset{R}}

\newcommand{\C}{\numberset{C}}

\journal{Computer Methods in Applied Mechanics and Engineering}

\begin{document}

\begin{frontmatter}
\title{Learning the Hodgkin-Huxley Model\\ with Operator Learning Techniques}

\author[label1]{Edoardo Centofanti}
\author[label1]{Massimiliano Ghiotto}
\author[label1]{Luca Franco Pavarino}

\affiliation[label1]{organization={Dipartimento di Matematica, Università di Pavia},
             addressline={Via Adolfo Ferrata, 5},
             city={Pavia},
             postcode={27100},
             country={Italy}
             }

\begin{abstract}
We construct and compare three operator learning architectures, DeepONet, Fourier Neural Operator, and Wavelet Neural Operator, 
in order to learn the operator mapping a time-dependent applied current to the transmembrane potential of the Hodgkin-Huxley ionic model.
The underlying non-linearity of the Hodgkin-Huxley dynamical system, the stiffness of its solutions, and the threshold dynamics depending on the intensity of the applied current, are some of the challenges to address when exploiting artificial neural networks  to learn this class of complex operators. 
By properly designing these operator learning techniques, we demonstrate their ability to effectively address these challenges, achieving  
a relative $L^2$ error as low as 1.4\% in learning the solutions of the Hodgkin-Huxley ionic model. 
\end{abstract}
\end{frontmatter}

\section{Introduction}
Scientific Machine Learning is a rapidly developing field combining the versatility of machine learning techniques, based on using data in order to learn or recognize patterns, with more traditional modelling techniques based on partial differential equations (PDE) \cite{ops21,gon23}. Scientific Machine Learning has been recently categorized into three main areas \cite{bou23}: PDE discovery, PDE solvers through artificial neural networks and Operator Learning.
The first area is related to inverse problems \cite{luk20}, while the last two are related to direct problems and they have been shown to be effective in reducing the solution time of many traditional PDE solvers.  Physics Informed Neural Networks (PINNs) and their variants \cite{cuo22,rai19,cai21,yan21,mah22} are popular examples of PDE solvers based on artificial neural networks: they have represented a revolution in the field, privileging a physics-oriented approach in optimizing the network parameters over a more traditional data-oriented approach. However, PINNs are trained in order to solve efficiently a single problem instead of a whole class of problems, where the solutions are the image of an operator depending on an input function over a given domain. In order to address this task, the third scientific machine learning branch, known as operator learning, has been developed over the last few years \cite{lu21,li20,tri23,mar23a,mar23b}. Operator learning has been applied to several problems, especially PDE systems depending on parameters or forcing terms, and also in industrial applications \cite{ors22, cao23}. Among the emerging approaches are DeepONet \cite{lu21} and neural operator \cite{kov21b} based architectures, such as Fourier \cite{li20} and Wavelet \cite{tri23} Neural Operators, which have shown remarkable capabilities in learning non-linear operators on regular or equispaced domains.

In this paper, we construct and study operator learning techniques to learn the dynamics of ionic models, a very important class of dynamical systems 
that describe the rich dynamics of excitable cells and biochemical reactions. These  models have significant applications in fields such as  Neuroscience to Computational Cardiology, see e.g. \cite{pav14,izh07,roc12,hod52,luo94,ten04}. The solutions of ionic models 
are strongly dependent on time-dependent currents applied to the system. In particular, they exhibits a threshold effect, with almost no response from the system when the intensity of the applied current is below a certain threshold, and a stiff nonlinear response for applied currents above threshold. 
Furthermore, the complex dynamics of ionic models can include the presence of limit cycles, bursting phenomena, and bifurcations.
These properties make this class of problems hard to learn by artificial neural networks, despite their well known non-linear approximation capabilities \cite{che95,dev21,hes97,she19}.

An interesting prototypical ionic model to study is undoubtedly the Hodgkin-Huxley model \cite{hod52}. It has been recently studied in the context of numerical machine learning in its fractional version with PINNs \cite{she23} and with multi-fidelity methods \cite{seu19}. Other recent works have recovered the Hodgkin-Huxley dynamics through recurrent maps built and trained as feed forward neural networks \cite{kup23}, through an extension of the neural ODE architecture \cite{che18, che20} and through a ResNet trained with short-term trajectories in the phase space \cite{su21}. Synchronization and oscillatory processes involving this model have also been studied through data driven machine learning approaches \cite{kem22, sag23}.    

In this work, we exploit operator learning techniques in order to learn the underlying non-linear operator, mapping time-dependent square-shaped applied currents 
to the Hodgkin-Huxley transmembrane potential. 
In particular, we construct three architectures, namely DeepONet, Fourier Neural Operator, and Wavelet Neural Operator, and study their performance in learning the Hodgkin-Huxley transmembrane potential dynamics.

This paper is organized as follows. In Section \ref{sec:HH}, we review the Hodgkin-Huxley model, its dynamics and the response of the system when different currents are applied. In Section \ref{sec:OL}, we present some mathematical details of operator learning and the architectures employed in this work and we briefly review the link between operator learning and differential equations theory. Finally, in Section \ref{sec:results}, we present our results and compare the learning strategies used, studying in particular how the generalization error changes with different architectures and their hyperparameters. 

\section{The Hodgkin-Huxley model}\label{sec:HH}
\noindent The Hodgkin-Huxeley (HH) model \cite{hod52, izh07} has been introduced in order to model the propagation of electrical signals through the squid giant axon and it is among the most studied biological dynamical systems in neuroscience and computational biology \cite{wan07, guc02, mcc07}. It has been recently source of inspiration for ANNs with neumorphic activations, such as the spiking neural networks \cite{maa97, kah22}.

Electrical activity in excitable cells (e.g. neurons or cardiac cells) is sustained and propagated via ionic currents through the cellular membrane. Different concentration of ionic species between the inside and the outside of the cell are responsible for electrochemical gradients, namely an electric potential between the intra and the extracellular space. The HH model starts from the hypothesis that the ionic currents are responsible for the generation of this potential through ionic-specific channels (also known as gates). The resulting currents across the membrane are strong enough to cause membrane depolarization. In particular, Hodgkin and Huxley determined experimentally that the current contributions are mainly three: a voltage-gated persistent $\chem{K}^+$ current with four activation gates, a voltage-gated transient $\chem{Na}^+$ current with three activation gates and one inactivation gate, and an Ohmic leak current, $I_L$, mostly carried by $\chem{Cl}^-$ ions. During the depolarization phase the membrane ionic permeability, related to sodium and potassium conductances, changes in order to restore the membrane potential. This re-polarizes the membrane through a sudden increase in the membrane potential (repolarization phase). These changes in the membrane permeability are due to the activation and inactivation of the ionic channels at a rate given by the so called \textit{gating variables}.
	
\noindent The 0D version of this model, describes the variation in time of the membrane voltage $V$ (in $mV$) through the axon membrane, when a current stimulus, $I_{app}$ (in $\mu A/cm^2$), is applied to the system. $V$ also depends on three gating variables, $m$, $h$ and $n$, which act as the probability for each ionic species to pass through the cellular membrane. 
	
\noindent Using the modern conventions employed in \cite{cle16}, the model can be written as the following ODE system:
	
	\begin{equation}\label{eq:HH}
		\begin{cases}
			C_m \frac{dV}{dt} = -\left( g_{Na} m^3h(V-V_{Na}) + g_Kn^4(V-V_K) + g_L(V-V_L) \right) + I_{app} \\[6pt]
			\displaystyle \frac{dm}{dt} = \alpha_m(V)(1-m)-\beta_m(V)m \\[8pt]
			\displaystyle \frac{dh}{dt} = \alpha_h(V)(1-h)-\beta_h(V)h \\[8pt]
			\displaystyle \frac{dn}{dt} = \alpha_n(V)(1-n)-\beta_n(V)n \\
		\end{cases}
	\end{equation}
	where we consider the following parameters and functions:
	\begin{align*}
		\alpha_m(V) &= 0.1(-V-50)\left[\exp\left(\frac{-V-50}{10}\right) - 1\right]^{-1}, & \!\!\beta_m(V) &\!\!= 4\exp\left(\frac{-V-75}{18}\right)\\
		\alpha_h(V) &= 0.07\exp\left(\frac{-V-75}{20}\right), &\!\!\!\!\beta_h(V) &\!\!= \!\!\!\left[ \exp\left( \frac{-V-45}{10} \right)\!+\!1 \right]^{-1} \\
		\alpha_n(V) &\!\!= 0.01(-V-65)\left[ \exp{\left(\frac{-V-65}{10}\right)}-1 \right]^{-1}, &\!\!\beta_n(V) &\!\!= 0.125\exp\left( \frac{-V-75}{80} \right)
	\end{align*}
	\begin{align*}
		&C_m = 1\,\mu F/cm^2,\,\, V_{Na} = 40\, mV,\,\, V_K = -87\, mV, \\
		&V_L = -64.387\, mV,\,\,  g_{Na} = 120\,mS/cm^2,\,\, g_K = 36\,mS/cm^2, \\
		&g_L = 0.3\, mS/cm^2,  
	\end{align*}
    \noindent
	and the initial conditions: $V_0 = -75 mV$, $n_0 = 0.317$, $m_0 = 0.05$ and $h_0 = 0.595$.

\noindent Various current inputs $I_{app}$ can be applied, resulting in different responses from the system. In our setup, we considered $I_{app}$ sampled from the set of single square functions with different intensity and width. An interesting property to point out is the translation invariance of the response: two pulses with the same intensity and width, but different starting time, will result in an identical solution up to a translation.

\begin{figure}
\centering
\includegraphics[width=\textwidth]{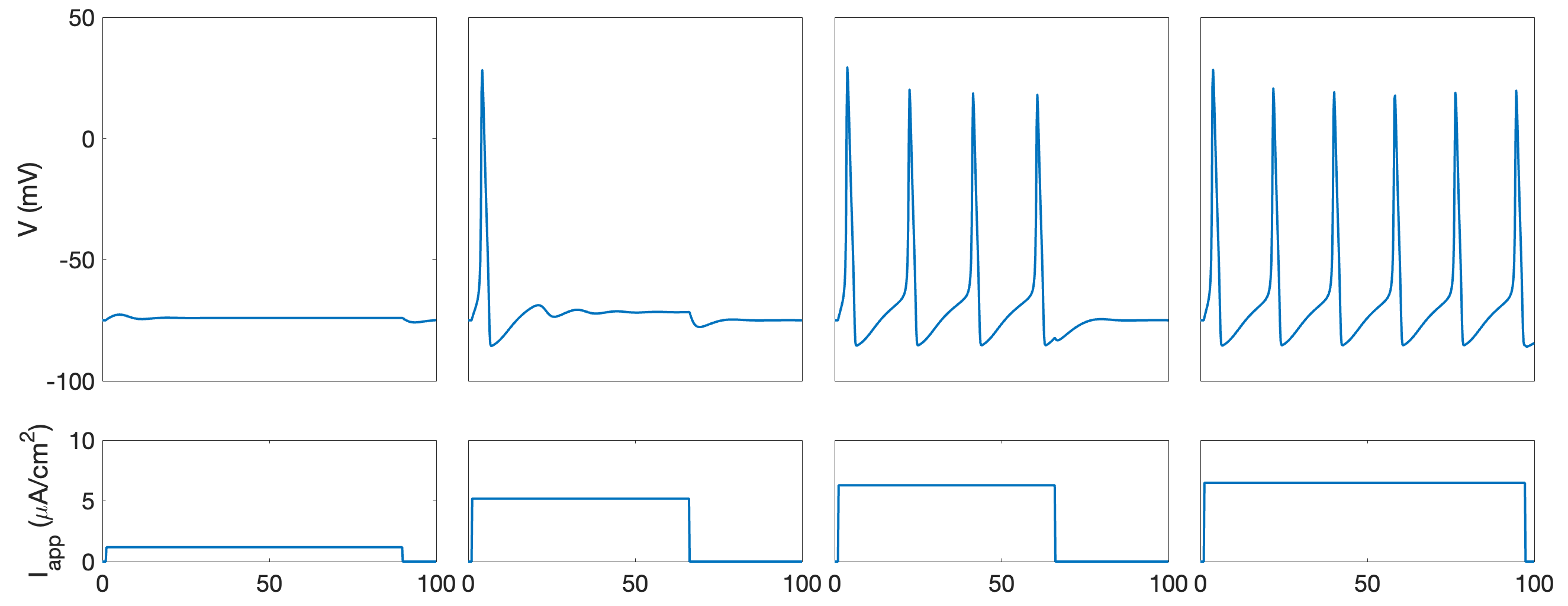}
\caption{Different $I_{app}$ pulses lead to different responses from the system. A pulse with an under-threshold height leads to zero peaks and a negligible perturbation of the system from its equilibrium state. Once the threshold value is crossed (around $2.21\, \mu A/cm^2 $ with these parameters), the system responds with a single peak or more peaks depending on the input current.}
\label{fig:pulses}
\end{figure}

\noindent We notice that the system (\ref{eq:HH}) can be seen as an integral operator mapping the time dependent applied current to the membrane potential $V$ \cite{kis97} and exploits the Volterra expansion for the reduced solution of an ODE system, which, in the case of the Hodgkin-Huxley model, reduces the 4 equations ODE system to an equivalent Integrate-And-Fire model for the membrane potential. Using the notation introduced in this work:

\begin{equation}\label{eq:volterra}
    V(t) = \eta(t-t^f) + \int_0^\infty \varepsilon^{(1)} (t-t^f;s)I_{app}(t-s)\,ds + \mathcal{O}(t^2)
\end{equation}
\noindent
where the expansion is stopped after the first order term, giving a good approximation for the solution $V(t)$ in (\ref{eq:HH}).
The kernels $\eta$ and $\varepsilon^{(1)}$ describe the form of the spike and the postsynaptic potential, respectively. Higher order terms depend on other kernels $\varepsilon^{(i)}$, where $i$ is the order, and describe the nonlinear interaction between several input pulses, but they are often negligible. We notice that the information related to the gating variables $m$, $h$ and $n$ is recovered through the kernels $\eta$ and $\varepsilon^{i}$, which are evaluated from the eigenvalues of the linearization of the system (\ref{eq:HH}) we started with.

\section{Operator Learning}\label{sec:OL}
\noindent Operator learning aims to approximate an unknown operator $G$, which often takes the form of the solution operator associated with a differential equation \cite{bou23}. Given pairs of data $(f,u)$, where $f\in \mathcal{A}$ and $u\in\mathcal{U}$ are from function spaces on domains $\Omega\in\mathbb{R}^d$ and $\Omega'\in\mathbb{R}^{d'}$ respectively, the goal is learning by mean of artificial neural networks, a potentially nonlinear and not always known operator $G:\mathcal{A}\rightarrow\mathcal{U}$ which maps $f$ to $u$. 

\noindent What is learned is an approximation $\hat{G}_\theta$ of the ground truth operator $G$. The approximation due to the neural network architecture employed, the optimizer used to minimize the network parameters and the number and location of the collocation points for the elements in the dataset are all potential sources of error \cite{mis22}.

\noindent Throughout this work, we will explore and compare the performance of operator learning with three architectures, namely Deep Operator Networks (DeepONets) \cite{lu21}, Fourier Neural Operators (FNOs) \cite{li20} and Wavelet Neural Operators (WNOs) \cite{tri23}.
	
\subsection{Deep Operator Networks}
\noindent The first theoretical result regarding the approximation property of nonlinear operators defined on infinite-dimensional spaces by shallow neural networks was provided in \cite{che95}. This result has been recently extended to deep neural networks (DNN) and a novel framework, called \textit{Deep Operator Network} (DeepONet), has been proposed in order to learn this class of operators \cite{lu21}. We report the following result:

\begin{teo}[Universal Approximation Theorem for Operator]
    Let $\sigma$ be a continuous non-polynomial function, $\mathcal{X}$ is a Banach space, $\mathcal{K}_1\subset \mathcal{X}$, $\Omega\subset \R^{d}$ are two compact sets in $\mathcal{X}$ and $\R^d$ respectively, $\mathcal{A}$ is a compact set in $C(\mathcal{K}_1)$, G a nonlinear continuous operator, mapping $\mathcal{A}$ into $\mathcal{U}\equiv C(\Omega)$. Then, for any $\varepsilon > 0$, there exist positive integers $n$,$p$ and $m$, and constants $c_i^k$, $\xi_{ij}^k$, $\theta_i^k$, $\zeta_k$ $\in \R$, $w_k\in \R^d$, $x_j\in \mathcal{K}_1$, $i=1,\cdots,n$, $k=1,\cdots,p$, and $j=1,\cdots,m$, such that
    \begin{equation}\label{eq:don_arch}
        \left| G(f)(y) - \sum_{k=1}^{p}\sum_{i=1}^{n}c_i^k \sigma \left( \sum_{j=1}^{m} \xi_{ij}^k f(x_j) + \theta_i^k \right) \sigma(w_k \cdot y + \zeta_k) \right| < \varepsilon
    \end{equation}
    holds for all $f\in \mathcal{A}$ and $y\in \Omega$. In this context, $C(\mathcal{K})$ is the Banach space of continuous functions defined on $\mathcal{K}$ with norm $||f||_{C(\mathcal{K})} = \max_{x\in K}|f(x)|$.
\end{teo} 

\noindent From (\ref{eq:don_arch}) we can identify the structure of the DeepONet architecture. The operator $G$ is approximated by a scalar product between two artificial neural networks, called \textit{branch network} and \textit{trunk network}, with the same output dimension $p$:

\begin{equation}
    \hat{G}_\theta(f)(y) = \sum_{k=1}^p \underbrace{b_k(f(x_1),...,f(x_{m}))}_{{}\text{branch net}}\underbrace{t_k(y)}_{{}\text{trunk net}}
    \label{eq:DeepOnetArc}
\end{equation}

\noindent Since our goal is learning operators taking functions as the input, we have to represent the input functions to the branch discretely, so that network training and consequent validation can be performed. This implies defining an \textit{encoder} \cite{lan22}:

\begin{deff}[Encoder]
    Given a set of sensor points $x_j\in \mathcal{K}_1$, for $j=1,\cdots, m$ we define the linear mapping
    \[ \mathcal{E}:\mathcal{A}\rightarrow\mathbb{R}^{m},\qquad \mathcal{E}(f) = (f(x_1),\cdots, f(x_m)), \]
    as the encoder mapping. The encoder $\mathcal{E}$ is well-defined as one can evaluate continuous functions pointwise.
\end{deff}
 
\noindent The only condition required for the encoding is that the sensor locations $\left\lbrace x_1, \dots, x_{m}\right\rbrace$ are the same but not necessarily on a lattice for all input functions $ f$, while we do not enforce any constraints on the output locations $y$.
\noindent In our numerical experiments, we actually implement a different encoding strategy, which will be detailed in section \ref{sec:results}. 

\begin{figure}
    \centering
    \includegraphics[trim={1cm 3cm 1cm 3cm},clip,width=0.8\textwidth]{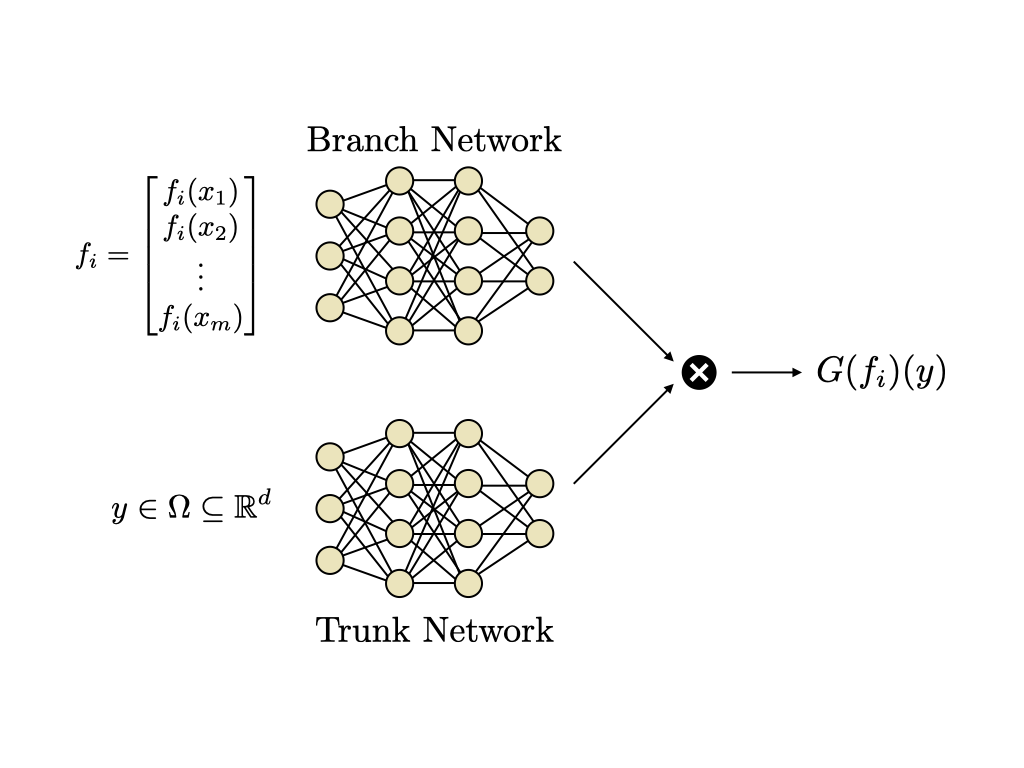}
    \caption{Schematic DeepONet architecture. The outputs of the branch and the trunk networks, here represented as forward neural networks, are multiplied through a scalar product in order to approximate the operator $G$.}
    \label{fig:don_arch}
\end{figure}
\noindent Finally, we define a \textit{reconstructor} (or \textit{decoder}) as the last operation performed in a DeepONet.
\begin{deff}[Reconstructor or Decoder]
    Given the trunk as a neural network $\tau: \R^n\rightarrow\R^{p+1}$, $y=(y_1,\cdots,y_n)\mapsto \{\tau_k(y)\}^p_{k=0}$, with $y\in U\subset \R^n$, we define a $\tau$-induced reconstructor (or decoder) as the map
    \begin{equation}
        D = D_\tau:\R^p\rightarrow C(U),\qquad D_\tau(\alpha_k):= \tau_0(y) + \sum_{k=1}^p\alpha_k\tau_k(y).
    \end{equation}
\end{deff}
\noindent The definition above formalizes the scalar product in (\ref{eq:DeepOnetArc}). We notice that each $\tau_k$ is the composition of linear operators and activation functions. $\tau_0$ is an extra bias often added in practical implementations, which in (\ref{eq:don_arch}) has been set to zero without loss of generality. 

\noindent The architecture introduced in \cite{lu21} is shown in Fig. \ref{fig:don_arch}. 
In their classical formulation, DeepONets are composed of DNN, but also other architectures involving Convolutional Neural Networks (CNNs), Residual Neural Networks (ResNets), and Recurrent Neural Networks (RNNs) like Gated Recurrent Unit (GRU) and Long Short Time Memory (LSTM) models are possible \cite{lu21,cao23,hex24}.
\subsection{Neural Operator}
\noindent The remaining two architectures which represent the object of the present work are based on the Neural Operator architecture \cite{kov21b}. Again, in the following we use $C(\Omega_t\times \Omega_t,\, \R^{d_{v_{t+1}}\times d_{v_{t}}})$ to denote the Banach space of continuous functions from $ \Omega_t \times \Omega_t$ to $\R^{d_{v_{t+1}}\times d_{v_{t}}}$ equipped with the sup-norm $||g||_\infty := \sup \{ |g(\tau)|:\tau\in\Omega_t\times\Omega_t \},$ $g\in C(\Omega_t\times \Omega_t,\, \R^{d_{v_{t+1}}\times d_{v_{t}}})$.
\begin{deff}[Neural Operator]
    Assuming that the input functions $f\in\mathcal{A}$ are $\R^{d}$-valued and defined on the bounded domain $\Omega\subset\R^d$ while the output functions $u\in\mathcal{U}$ are $\R^{d_u}$ valued and defined on the bounded domain $\Omega'\subset\R^{d'}$, we define the Neural Operator as the architecture $\hat{G}_\theta:\mathcal{A}\rightarrow\mathcal{U}$ with the following structure:
    \begin{itemize}
        \item \textbf{Lifting:} Through a pointwise function $\R^{d}\rightarrow\R^{v_0}$, map the input $\{f: \Omega\rightarrow \R^{d}\}\mapsto\{v_0: \Omega\rightarrow\R^{d_{v_0}}\}$ to its first hidden representation. A suitable choice is $d_{v_0} > d$, thus this is called lifting operation performed by a fully local operator.
        \item \textbf{Iterative Kernel Integration:} For $t=0,\cdots,T-1$, each hidden representation is mapped to the next, namely $\{v_t:\Omega_t\rightarrow \R^{d_{v_t}}\}\mapsto\{v_{t+1}:\Omega_{t+1}\rightarrow\R^{d_{v_{t+1}}}\}$, through the action of the sum of a local linear operator, a non-local integral kernel operator, and a bias function, composing the sum with a fixed pointwise nonlinearity. Here we set $\Omega_0=\Omega$ and $\Omega_T=\Omega'$ and impose that $\Omega_t\subset\R^{d_t}$ is a bounded domain.
        \item \textbf{Projection:} Using a pointwise function $\R^{d_{v_{T}}}\rightarrow\R^{d_u}$, map the last hidden representation $\{v_T:\Omega'\rightarrow\R^{d_{v_{T}}}\}\mapsto\{u:\Omega'\rightarrow\R^{d_u}\}$ to the output function. In this case the suitable choice of dimensionality is $d_{v_{t}}>d_u$, thus this is a projection step performed by a fully local operator.
    \end{itemize}
\end{deff}
\noindent From the definition above, we thus have
\begin{equation}\label{eq:neural_op}
    \hat{G}_\theta := Q\circ \sigma_T\left( W_{T-1}+K_{T-1}+b_{T-1} \right)\circ \cdots\circ \sigma_1\left( W_0+K_0 + b_0 \right)\circ P
\end{equation}
with $P:\R^{d}\rightarrow\R^{d_{v_{0}}}$,  $Q:\R^{d_{v_{T}}}\rightarrow\R^{d_u}$ are the local lifting and projection operators, $W_t\in\R^{d_{v_{t+1}}\times d_{v_{t}}}$ are local linear operators (matrices), $K_t : \{v_t:\Omega_t\rightarrow \R^{d_{v_{t}}}\} \rightarrow \{v_{t+1}:\Omega_{t+1}\rightarrow \R^{d_{v_{t+1}}}\}$ are integral kernel operators, $b_t:\Omega_{t+1}\rightarrow\R^{d_{v_{t+1}}}$ are bias functions and $\sigma_t$ are fixed activation functions acting locally as maps $\R^{v_{t+1}}\rightarrow\R^{v_{t+1}}$. Input and output dimensions $d_{v_{0}},\cdots,d_{v_{T}}$ as well as the domains of definition $\Omega_1,\cdots, \Omega_{T-1}$ are hyperparameters of the architecture. Furthermore, as local maps we consider their action as pointwise, namely for lifting and projection we have respectively $(P(a))(x) = P(a(x))$ for $x\in D$ and $(Q(a))(x) = Q(a(x))$ for $x\in \Omega'$. For the activation, we have instead $(\sigma_t(v_{t+1}))(x) = \sigma(v_{t+1}(x))$ for any $x\in \Omega_{t+1}$.

\noindent We now provide three admissible definitions for the integral kernel operators employed in (\ref{eq:neural_op}).
\begin{deff}[Integral Kernel Operators]
    The following definitions for the integral kernel operators $K_t$ are admitted \cite{kov21b}:
    \begin{itemize}
        \item Let $\kappa^{(t)}\in$ $C(\Omega_{t+1}\times \Omega_t; \R^{d_{v_{t+1}}\times d_{v_{t}}})$ and let $\nu_t$ be a Borel measure on $\Omega_t$ (normally simply the Lebesgue measure on $\R^{d_t}$). $K_t$ is defined by
        \begin{equation}\label{eq:kappat}
            (K_t(v_t))(x) = \int_{\Omega_t} \kappa^{(t)}(x,y)v_t(y)\,d\nu_t(y),\qquad \forall x \in \Omega_{t+1}.
        \end{equation}
        \item Let $\kappa^{(t)}\in$ $C(\Omega_{t+1}\times \Omega_t\times \R^{d}\times\R^{d}; \R^{d_{v_{t+1}}\times d_{v_t}})$. $K_t$ is defined by
        \begin{equation}
             (K_t(v_t))(x) = \int_{\Omega_t} \kappa^{(t)}(x,y,f\left(\Pi_{t+1}(x)\right),f\left(\Pi_{t}(y)\right))v_t(y)\,d\nu_t(y),\qquad \forall x \in \Omega_{t+1}.
        \end{equation}
        where $\Pi_t:\Omega_t\rightarrow \Omega_t$ are fixed projection maps.
        \item Let $\kappa^{(t)}\in$ $C(\Omega_{t+1}\times \Omega_t\times \R^{d_{v_t}}\times\R^{d_{v_t}}; \R^{d_{v_{t+1}}\times d_{v_t}})$. $K_t$ is defined by
        \begin{equation}
            (K_t(v_t))(x) = \int_{\Omega_t} \kappa^{(t)}(x,y,v_t\left(\Pi_{t}(x)\right),v_t\left(y\right))v_t(y)\,d\nu_t(y),\qquad \forall x \in \Omega_{t+1}.
        \end{equation}
        where $\Pi_t:\Omega_t\rightarrow \Omega_t$ are fixed projection maps, often taken to be the identity.
    \end{itemize}
\end{deff}
\noindent From now on, we will consider $K_t$ as in the first definition. We now need to define the single hidden layer update
\begin{deff}[Single Hidden Layer Update]
    Given $K_t$ as in (\ref{eq:kappat}), we define the single hidden layer update rule as
    \begin{equation}
        v_{t+1}(x) = \sigma_{t+1} \left( W_tv_t\left( \Pi_t(x) \right) + \int_{\Omega_t} \kappa^{(t)}(x,y)v_t(y)\,d\nu_t(y) + b_t(x)\right)\qquad \forall x\in \Omega_{t+1}.
    \end{equation}
    where $\Pi_t:\Omega_{t+1}\rightarrow \Omega_t$ are fixed projection maps.
\end{deff}

\noindent Having defined the general structure for the Neural Operator, we can now describe the FNO and the WNO by simply specializing the integral kernel (\ref{eq:kappat}) for each architecture.

\subsubsection{Fourier Neural Operator}
\begin{figure}
    \centering
    \includegraphics[width=0.8\textwidth]{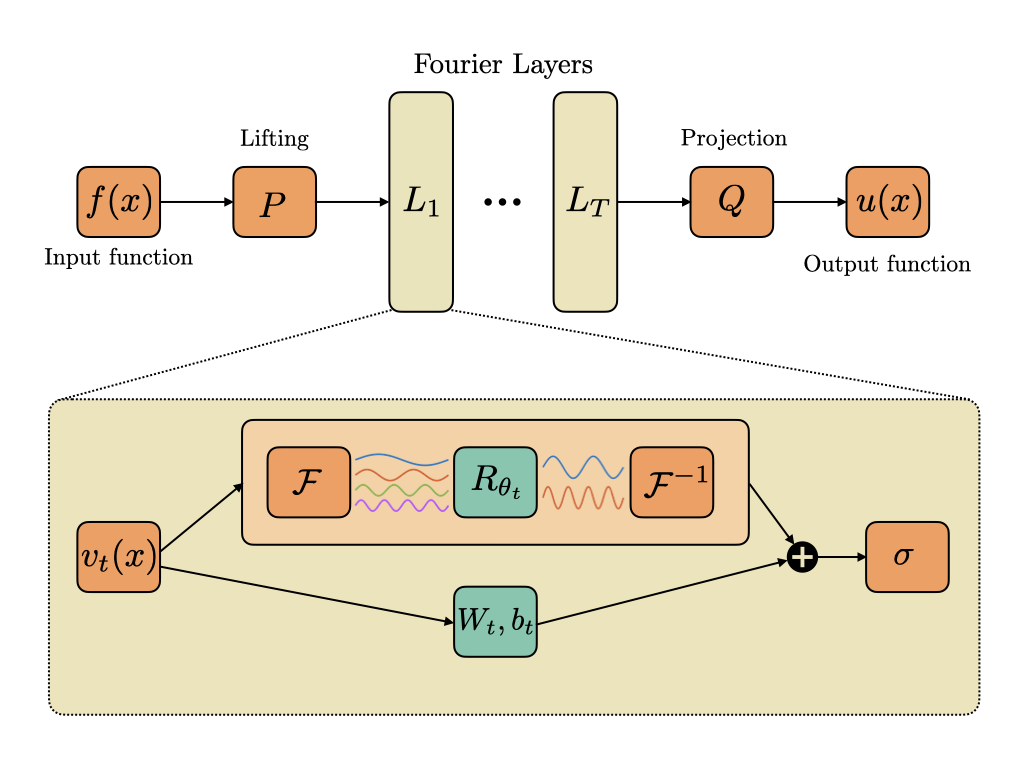}
    \caption{Schematic FNO architecture.}
    \label{fig:fno_arch}
\end{figure}

\noindent
The FNO works with the integral kernel parameterized in the Fourier space. We thus need to define the Fourier transform and anti-transform.
\begin{deff}[Fourier transform and anti-transform]
    Let $\mathcal{F}:L^2(\Omega,\mathbb{C}^n)\to \ell^2(\mathbb{Z}^d,\mathbb{C}^n)$ denote the Fourier transform and $\mathcal{F}^{-1}$ its inverse. Given $u\in L^2(\Omega,\mathbb{C}^n)$ and $\hat{u}\in\ell^2(\mathbb{Z}^d,\mathbb{C}^n)$, we define
\begin{equation}
    \begin{split}
        & \left( \mathcal{F}u \right)_j(k) =\hat{u}_j(k) = \langle u_j, \varphi_k\rangle_{L^2(\Omega,\mathbb{C}^n)},\qquad j\in\{1,\cdots,n\},\quad k\in\mathbb{Z}^d, \\
        &\left( \mathcal{F}^{-1} \hat{u} \right)_j(x) = \sum_{k\in\mathbb{Z}^d}\hat{u}_j(k)\varphi^{-1}_k(x),\qquad j\in\{1,\cdots,n\},\quad x\in\Omega \\
        &\varphi_k(x) := e^{-2\pi ik \cdot x},\qquad x\in \Omega
    \end{split}
\end{equation}
\end{deff}
\noindent
For the FNO, $\Omega_t = \mathbb{T}^d$, namely the domains considered for each layer are the unit torus. We also assume that the integral operator is in the form (\ref{eq:kappat}), i.e. it does not depend on the input function $f$ or on the output of each single hidden layer update $v_t$. In particular, let us consider $\kappa^{(t)} \in C(\Omega_t\times \Omega_t,\, \R^{d_{v_{t+1}}\times d_{v_{t}}})$ such that
\begin{equation*}
(K_t(v_t))(x) = \int_{\Omega_t} \kappa^{(t)}(x, y) v_t(y) \ dy \quad \forall x \in \Omega.
\end{equation*}
\noindent
Here, the integral kernel $\kappa^{(t)} : \Omega_t\times\Omega_t \to \R^{d_{v_{t+1}} \times d_{v_{t}}}$ is a function parameterized by some parameters $\theta_t$ in a suitable parameters space $\Theta_{t}\subset\R^{d_{v_{t+1}}\times d_{v_{t}}}$. We can thus also identify $\kappa^{(t)}=\kappa^{(t)}_{\theta_t}$. Supposing that $\kappa^{(t)}_{\theta_t}(x, y) = \kappa^{(t)}_{\theta_t}(x-y)$, we can apply the convolution theorem and write
\begin{equation}\label{eq:FNOKernel}
    (K_t(v_t))(x) = \int_{\Omega} \kappa^{(t)}_{\theta_t}(x-y)v_t(y)\,dy = \mathcal{F}^{-1}\left( \mathcal{F}\left(\kappa^{(t)}_{\theta_t}\right)\cdot\mathcal{F}(v_t) \right)(x)\qquad \forall x\in\Omega_t
\end{equation}
\noindent
Where the dot operation is defined as
\begin{equation} \label{eq:hadamard}
    \mathcal{F}\left(\kappa^{(t)}_{\theta_t}\right) \cdot\mathcal{F}(v_t) := \sum_{k \in \Z} \mathcal{F}\left(\kappa^{(t)}_{\theta_t}\right)(k)\cdot \mathcal{F}(v_t)(k) = \sum_{k \in \Z} \hat{\kappa}^{(t)}_{\theta_t}(k) \cdot \hat{v}_t(k).
\end{equation}
\noindent Taken $k\in\mathbb{Z}^d$ as frequency mode, for $k$ fixed we have that $\mathcal{F}\left(\kappa^{(t)}_{\theta_t}\right)(k)\in\mathbb{C}^{d_{v_{t+1}}\times d_{v_{t}}}$ and $\mathcal{F}(v_t)(k)\in\mathbb{C}^{d_{v_{t}}}$. This means that $\kappa^{(t)}_{\theta_t}$ can be parameterized directly by its Fourier coefficients, thus,
\begin{equation}\label{eq:KernelFNO}
   (K_{t}(v_t))(x) = \mathcal{F}^{-1}\left( R_{\theta_t}\cdot\mathcal{F}(v_t) \right)(x),\qquad\forall x \in \Omega_t.
\end{equation} 
where $R_{\theta_t} = \hat{\kappa}^{(t)}_{\theta_t}(k)$ for $k$ fixed. Fig. \ref{fig:fno_arch} reports a scheme relative to the FNO architecture considered in this work.

\noindent We must make some approximations from continuous space to discrete ones to have a finite-dimensional parameterization. This approximation of the FNO is called pseudo$(\Psi)$-Fourier Neural Operator ($\Psi$-FNO) \cite{kov21a}. In this case, the domain $\Omega_t$ is discretized with $J\in\mathbb{N}$ points, thus $v_t$ can be treated as a tensor in $\mathbb{C}^{J\times d_{v_{t}}}$. Since integrals cannot be calculated exactly, this leads to considering the Fourier series truncated at a maximal mode $k_{\max}$, such that
\begin{equation*}    
    k_{\max} = |Z_{k_{\max}}| = \big|\{k\in\mathbb{Z}^d \ :\ |k_j|\le k_{\max,j},\text{ for } j=1,\dots,d\}\big|.
\end{equation*}
\noindent
In practical implementation, the Fourier transform is replaced by the Fast Fourier Transform (FFT), and the weight tensor $R_{\theta_t}$ is parameterized as a complex-valued tensor $ R_{\theta_t} \in \C^{k_{\max}\times d_{v_{t+1}}\times d_{v_{t}}}$ for $i = 1, \dots, L$. From now on we represent $R_{\theta_t}$ with only $R$ to simplify the notation. 
Finally, we note that if $v_t$ is real-valued and we want that $v_{t+1}$ to be real-valued too, we can enforce conjugate symmetry in tensor $R$, namely $R(-k)_{j,l}=\overline{R}(k)_{j,l}$. 

\subsubsection{Wavelet Neural Operator}
\noindent
The Wavelet Neural Operator \cite{tri23} learns the network parameters in the wavelet space, which is both frequency and spatial localized. The WNO has a similar structure to FNO but exploits wavelet integrals instead of Fourier integrals.

\noindent Again, we impose translation invariance for $\kappa^{(t)}_{\theta_t}$ and we can rewrite our integral kernel (\ref{eq:kappat}) as

\begin{equation}
    (K_{t}(v_t))(x) = \int_{\Omega} \kappa^{(t)}_{\theta_t}(x-y)v_t(y)\,dy\qquad \forall x\in\Omega.
\end{equation}
\noindent As in the previous case, we exploit the convolution theorem, but we apply it in order to parameterize $\kappa^{(t)}_{\theta_t}$ in the wavelet domain.

\noindent A function $\psi\in L^2(\Omega,\R)$ is called \textit{orthonormal mother wavelet} if it can be used to define a complete orthonormal basis for $L^2(\Omega,\R)$. A set of orthonormal functions for this space in this setting is built from the mother wavelet through scaling and shifting operations, obtaining a family of functions $\{\psi_{k,l}: k,l\in\mathbb{Z}\}$ with the following form:
\begin{equation}\label{eq:waveletBasis}
    \psi_{k,l}(x) = 2^{\frac{k}{2}}\psi\left( 2^kx-l \right) 
\end{equation}

\begin{figure}
    \centering
    \includegraphics[trim={0cm 3cm 0cm 0cm},clip,width=0.8\textwidth]{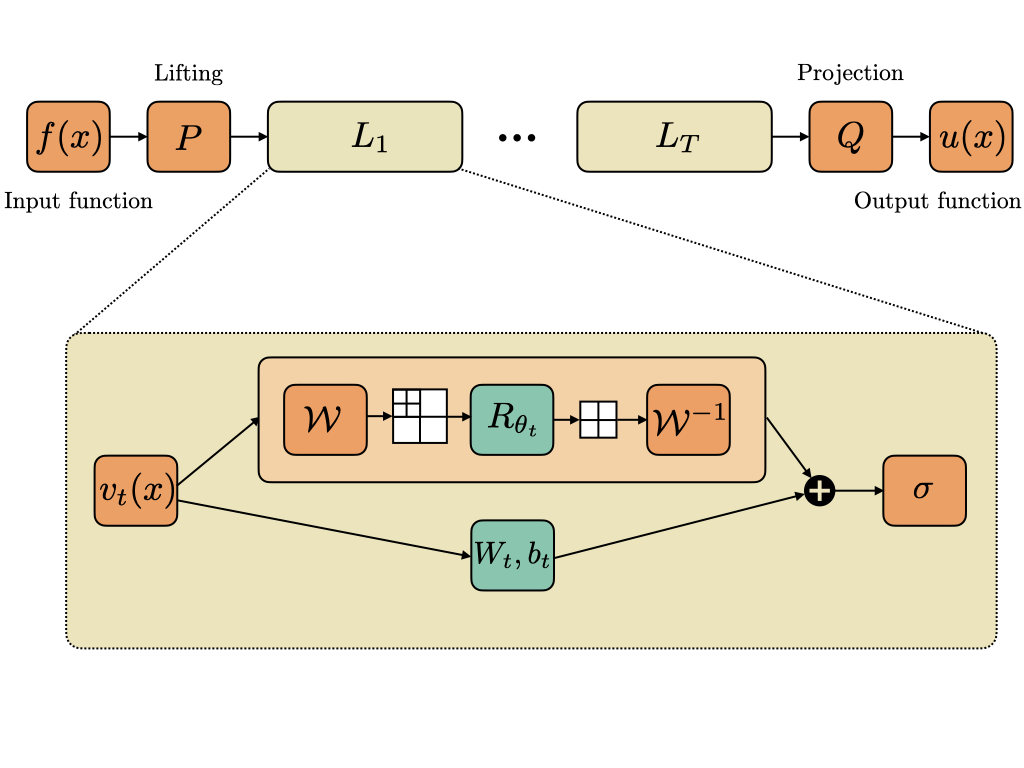}
    \caption{Schematic WNO architecture.}
    \label{fig:wno_arch}
\end{figure}
\noindent We can now define $\mathcal{W}$ and $\mathcal{W}^{-1}$, the forward and inverse wavelet transform.
\begin{deff}
    Given $g\in L^2(\Omega,\R^n)$, $h\in L^2(\R^+_*\times\R,\R^n)$, $\alpha\in\R^+_*$ (i.e. $\R^+$ including 0) and $\beta\in\R$, we define the forward wavelet transform $\mathcal{W}(g)$ as
    \begin{equation}
        (\mathcal{W}g)_j(\alpha,\beta) = \int_\Omega g_j(x) \frac{1}{|\alpha|^{\frac{1}{2}}}\psi\left(\frac{x-\beta}{\alpha}\right)d\,x
    \end{equation}
    and the inverse wavelet transform $\mathcal{W}^{-1}(h)$ as
    \begin{equation}\label{eq:CWT}
        (\mathcal{W}^{-1}h)_j(x) = \frac{1}{C_\psi}\int_0^\infty \int_\Omega h_j(\alpha,\beta) \frac{1}{|\alpha|^{\frac{1}{2}}}\tilde{\psi}\left( \frac{x-\beta}{\alpha} \right)d\beta\frac{d\alpha}{\alpha^2}
    \end{equation}
    where $\psi\in L^2(\Omega,\R)$ is the mother wavelet, $\tilde{\psi}$ is called dual mother wavelet (for our case $\psi=\tilde{\psi}$) and $C_\psi$ is defined as
    \begin{equation}
        C_\psi = 2\pi \int_\Omega \frac{|\psi(k)|^2}{|k|}\,dk
    \end{equation}
\end{deff}
\noindent In the previous definition, $\alpha$ and $\beta$ are scaling and translation parameters respectively. From (\ref{eq:waveletBasis}), $\alpha = 2^{-k}$ and $\beta = 2^{-k}l$, thus for each choice of the pair $(k,l)$ exists a pair $(\alpha,\beta)$ and we can map $\psi_{k,l}\mapsto\psi_{\alpha,\beta}:=\psi\left( \frac{x-\beta}{\alpha} \right)$. 

\noindent From the convolution theorem, one can prove also for the wavelet transform case, as in (\ref{eq:FNOKernel}), that
\begin{equation}
    K_t(x) = \mathcal{W}^{-1}\left( R_{\theta_t}\cdot \mathcal{W}(v_t) \right)(x);\quad x\in \Omega,
\end{equation}
where $R_{\theta_t}$ in this case is the wavelet transform of $\kappa^{(t)}_{\theta_t}$ parameterized by some in a suitable parameters space $\Theta_{t}\subset\R^{d_{v_{t+1}}\times d_{v_{t}}}$.

\noindent Analogously to the $\Psi$-FNO case, in order to implement the WNO, one has to consider the discrete wavelet transform (DWT). Moreover, $R_{\theta_t}$ is parameterized as a real-valued $(\zeta\times d_{v_{t+1}}\times d_{v_t})$ tensor, where $\zeta = n_s/2^k + (2n_w-2)$, with $n_s$ length of the discrete input of the DWT and $n_w$ length of the vanishing moment of the undertaken wavelet, after having applied recursively $k$ times the DWT. In practice, we fix the level $k$ in (\ref{eq:waveletBasis}) such that the corresponding $\psi_{k,l}$ contain the most important features of the transformed function. This usually happens when $\psi_{k,l}$ is related to the low frequency information, corresponding to coefficients of (\ref{eq:CWT}) with higher scales $2^k$. Figure \ref{fig:wno_arch} shows schematically the WNO architecture described in this work. For more information about signal processing and wavelet transform we suggest \cite{del95,mal99,fra06}.

\section{Results} \label{sec:results}
\noindent In this section we present the results of our numerical tests. We have considered the architectures introduced in Section \ref{sec:OL}, and ran experiments changing the hyperparameters setup, the choices of activation function and the network complexity, namely the number of neurons and trainable parameters. All the tests presented in this work have been performed exploiting the \texttt{PyTorch} library running on an Nvidia Quadro RTX $5000$ GPU, with $16$ GB of memory.

\subsection{Dataset, loss and optimizer} \label{sec:dataset}
\noindent For our dataset, we have considered 2000 numerical solutions of the problem (\ref{eq:HH}), divided between 1600 training data and 400 test data. During the training phase the dataset has been divided into mini batches of size 100. The numerical solutions in the dataset have been obtained from the same number of $I_{app}$ applied to the system, randomly chosen. In particular, the square pulses have been extracted with a uniformely distributed duration over a timespan of 100 $ms$ and with a uniformely distributed height from 0 to 10 $\mu A/cm^2$. The solutions have been sampled on an equispaced 500 points grid and solved numerically using the MATLAB routine \texttt{ode15s}.

\noindent Since our input functions are square pulses in time, the encoder for DeepONet can be reconstructed in the following way:
\begin{equation}
    f(x) = \hat{I}\,H(x-t_i)\,H(t_f-x)
\end{equation}
where $H$ is the Heaviside step function, $\hat{I}$ is the intensity of the pulse and $t_i$ and $t_f$ are the starting and the ending time of the pulse. Considering as encoding strategy the function $\hat{\mathcal{E}}:\mathcal{A}\rightarrow\R^3$ defined as
\begin{equation}
    \hat{\mathcal{E}}(f) = (t_i,t_f,\hat{I})
\end{equation}
it is possible to save memory resources without losing information, since the triple $(t_i,t_f,\hat{I})$ completely characterizes a single input function. For FNO and WNO we have instead considered as input the whole signal discretized in time. As a decoder for the DeepONet we have considered the dot product between the trunk and the branch network outputs. As loss for all the architectures explored, we have chosen the relative $L^2$ error, namely
\begin{equation}\label{eq:rel2}
    L_2 = \frac{1}{N}\sum_{i=1}^N \frac{|| \hat{G}_\theta(f_i)-u_i ||_{L^2(\Omega)}}{||u_i||_{L^2(\Omega)}}.
\end{equation}
\noindent
where $\hat{G}_\theta$ is the operator approximated by the architecture of our choice, while $f_i, u_i$ are respectively the $i$-th input function and solution of the dataset.

\noindent Regarding the parameters optimization, we considered the PyTorch \texttt{ReduceOnPlateau} policy as learning rate strategy, where the starting learning rate ($1\cdot 10^{-3}$ for DeepONet and WNO, $7\cdot10^{-4}$ for FNO) is reduced of 10\% when the loss is non-decreasing for 10 epochs and we have considered the AdamW algorithm as optimizer.

\subsection{DeepONet results}
\noindent In all the DeepONet tests performed, as input for the branch, we have considered the output of the encoder discussed in the previous section, while for the trunk, we have considered the 500 equispaced domain points.

\noindent For both branch and trunk networks, we have considered feed-forward neural networks with layer normalization \cite{ba16} and adaptive activation function \cite{jag20}. Both branch and trunk networks employ the \texttt{tanh} activation function. 

\noindent In Table \ref{table:results} we report our best results for this architecture. The best mean test loss, evaluated from three trained networks each one with different randomly initialized parameters is around 2.2\% and has been obtained with a branch and trunk network width of 700 neurons, for a total of 3.4 million parameters.

\begin{figure}
\begin{subfigure}[b]{0.5\textwidth}
    \centering
    \includegraphics[width=\textwidth]{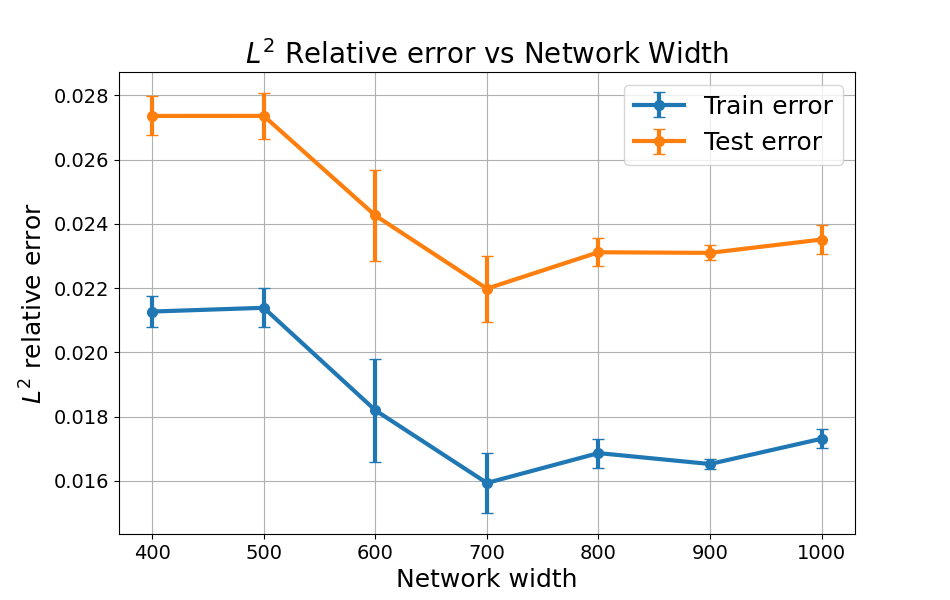}
    \caption{}
    \label{fig:don_scatter}
\end{subfigure}%
\begin{subfigure}[b]{0.5\textwidth}
    \centering
    \includegraphics[width=\textwidth]{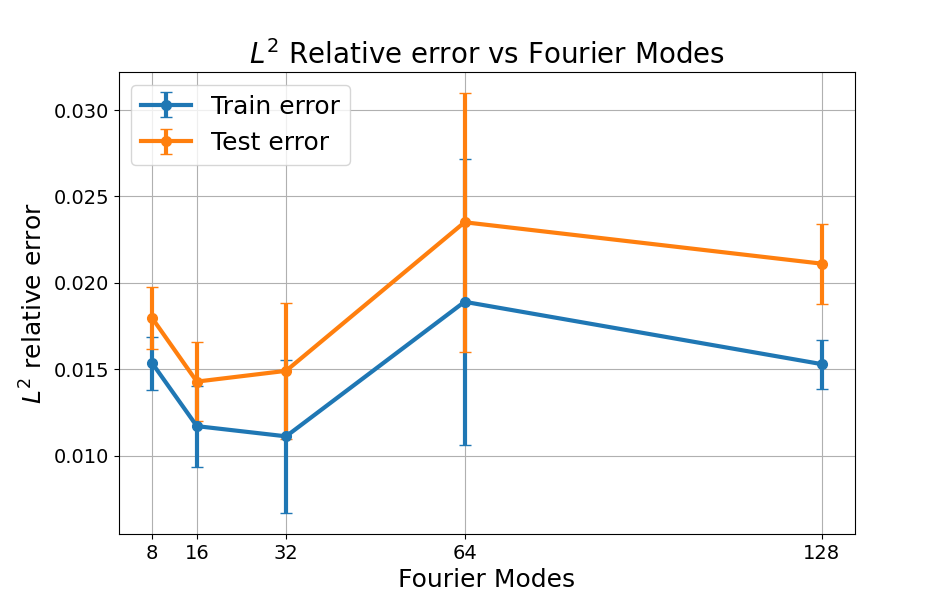}
    \caption{}
    \label{fig:fno_scatter}
\end{subfigure}
\centering
\begin{subfigure}[b]{0.5\textwidth}
    \includegraphics[width=\textwidth]{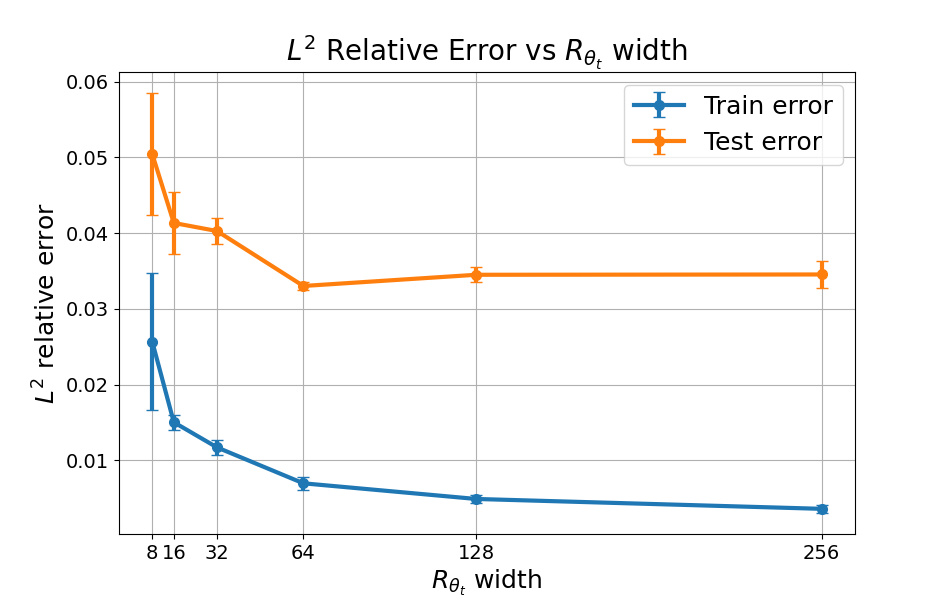}
    \caption{}
    \label{fig:wno_scatter}
\end{subfigure}
    \caption{ (\ref{fig:don_scatter}) $L^2$ train and test loss vs network width for DeepONet. All the branch nets considered have a depth of 4, while all the trunks considered have a depth of 3. The width is indicated on the $x$ axis and it is the same for all the layers and the output layer. A width of 700 neurons is the best among the ones explored for our problem. The training error is indicated with a blue line, while the test error is indicated with an orange line. (\ref{fig:fno_scatter}) $L^2$ train and test loss vs Fourier modes for the Fourier transform employed in the FNO. 16 and 32 modes resulted in a better performing trained network for our problem with our configuration. The training error is indicated with a blue line, while the test error is indicated with an orange line. (\ref{fig:wno_scatter}) $L^2$ train and test loss vs $R_{\theta_t}$ width for WNO. $R_{\theta_t}$ is a tensor of trainable parameters in the wavelet space. The width is indicated on the $x$ axis and it is the same for all the levels of the wavelet transform. In our numerical tests, a width of 64 has resulted in the best approximation error. The training error seems to improve with a superlinear rate when the width is increased. The training error is indicated with a blue line, while the test error is indicated with an orange line.}
\end{figure}

\noindent In Fig. \ref{fig:don_scatter} we study the behavior of the train and test losses as a function of the width of the networks, which is the same among all the layers, from 400 up to 1000 neurons. After an absolute minimum in the test error corresponding to a network width of 700 neurons, the error increases again up to 2.4\%. 

\noindent In Fig. \ref{fig:DON_test} we report four examples of solutions approximated by DeepONet, for our best model. The four corresponding inputs and the numerical solutions are taken from the test dataset. Single or double peak solutions are in general very well approximated, with the $L^2$ relative error below the mean test error. However, solutions with many peaks (e.g. the six-peaks solution shown in figure) tend to accumulate the error and the approximation deteriorates in time. 

\begin{table}
    \centering
    \begin{tblr}{
        colspec = {cccccc},
        }
        \toprule
        \SetCell[r=9]{c} \rotatebox[origin=c]{90}{DeepONet} & Test Error & \SetCell[c=2]{c} Configuration & & Total Parameters & Epochs \\
        \cmidrule{3-4}
         & & Trunk & Branch &  &  \\
        \midrule
        & 0.0274$\pm$0.0006 & [400]*3 & [400]*4   & 1.1M  & \SetCell[r=7]{c} $10000$  \\
        & 0.0274$\pm$0.0007 & [500]*3  & [500]*4   & 1.7M &  \\
        & 0.0243$\pm$0.0014 & [600]*3  & [600]*4   & 2.5M & \\
        & 0.0220$\pm$0.0010 & [700]*3  & [700]*4   & 3.4M & \\
        & 0.0231$\pm$0.0004 & [800]*3  & [800]*4   & 4.5M &   \\
        & 0.0231$\pm$0.0002 & [900]*3  & [900]*4   & 5.6M &  \\
        & 0.0235$\pm$0.0005 & [1000]*3 & [1000]*4 & 7.0M  &       \\
        \bottomrule
        \SetCell[r=6]{c} \rotatebox[origin=c]{90}{FNO} & & \SetCell[c=2]{c} Fourier modes & &  &   \\
        \midrule
        & 0.0180$\pm$0.0018 & \SetCell[c=2]{c} 8   & & 34k  & \SetCell[r=5]{c} $2000$ \\
        & 0.0143$\pm$0.0023 & \SetCell[c=2]{c} 16  & & 58k  &  \\
        & 0.0149$\pm$0.0039 & \SetCell[c=2]{c} 32  & & 108k & \\
        & 0.0235$\pm$0.0075 & \SetCell[c=2]{c} 64  & & 206k &  \\
        & 0.0211$\pm$0.0023 & \SetCell[c=2]{c} 128 & & 402k & \\
        \bottomrule
        \SetCell[r=7]{c} \rotatebox[origin=c]{90}{WNO}& & \SetCell[c=2]{c} $R_{\theta_t}$ width & &  & \\
        \midrule
        & 0.0505$\pm$0.0080 & \SetCell[c=2]{c} 8   & & 206k & \SetCell[r=6]{c} $2000$\\
        & 0.0413$\pm$0.0040 & \SetCell[c=2]{c} 16  & & 821k & \\
        & 0.0403$\pm$0.0018 & \SetCell[c=2]{c} 32  & & 3.2M &  \\
        & 0.0330$\pm$0.0005 & \SetCell[c=2]{c} 64  & & 13M  &  \\
        & 0.0345$\pm$0.0010 & \SetCell[c=2]{c} 128 & & 52M  &  \\
        & 0.0345$\pm$0.0018 & \SetCell[c=2]{c} 256 & & 209M & \\
        \bottomrule
        \end{tblr}
\caption{Our best results obtained from the implementation of DeepONet, FNO and WNO for the prediction of the membrane potential in the Hodgkin-Huxley model. For this problem the best relative validation error for DeepONet is around $2.1\%$, penalized by solutions resulting in more than 3 peaks, which have been particularly challenging to learn. The notation $[w]*L$ indicates a FNN composed by $L$ layer of width $w$. For FNO the best error is around $1.4\%$, obtained with 16 Fourier modes. For WNO this value is $3.3\%$ and the parameterized kernel has been discretized into a tensor $R_{\theta_t}$ with 64 as first dimension (width).} \label{table:results}
\end{table}

\subsection{FNO results}
\noindent In all the FNO tests performed, we considered all the $500$ discretization points as input to the network, without further encoding.

\noindent In Table \ref{table:results} we report our best results for this architecture. Among the hyper-parameters of the best performing architecture we highlight: lifting dimension $d_1=32$, Maximum number of Fourier modes considered for the FFT $k_{max}=16$, number of Fourier integral layers $L=3$, Rectified Linear Unit (ReLU) activation function. The best mean test loss, evaluated from three trained networks each one with different randomly initialized parameters is around 1.4\% and has been obtained considering 16 Fourier modes in the integral kernel, for a total of 58k parameters.

\noindent In Fig. \ref{fig:fno_scatter} we study how the relative $L^2$ error varies when we change the number of Fourier modes for the FFT employed in the approximation of the FNO integral kernel. After observing an initial decline, extending beyond 32 modes resulted in a rise in both test and training errors. 

\begin{figure}
    \centering
    \begin{subfigure}[b]{\textwidth}
        \centering
        \includegraphics[width=\textwidth]{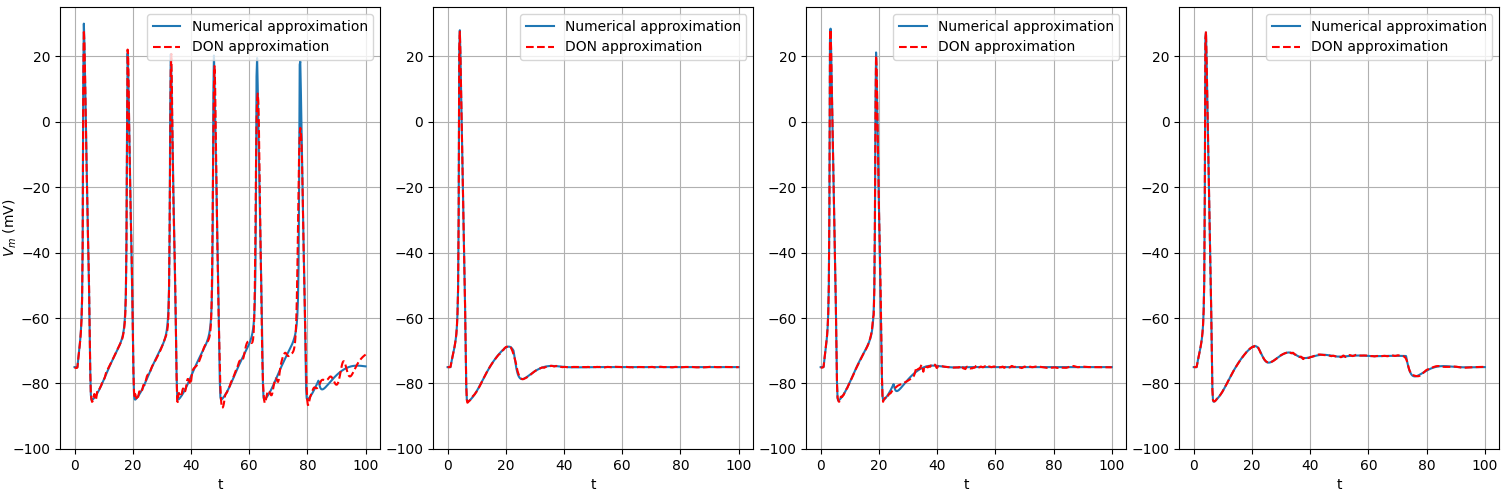}
        \caption{Deep Operator Network approximations.}
        \label{fig:DON_test}
    \end{subfigure}
    \vfill \vspace{0.35cm}
    \begin{subfigure}[b]{\textwidth}
        \centering
        \includegraphics[width=\textwidth]{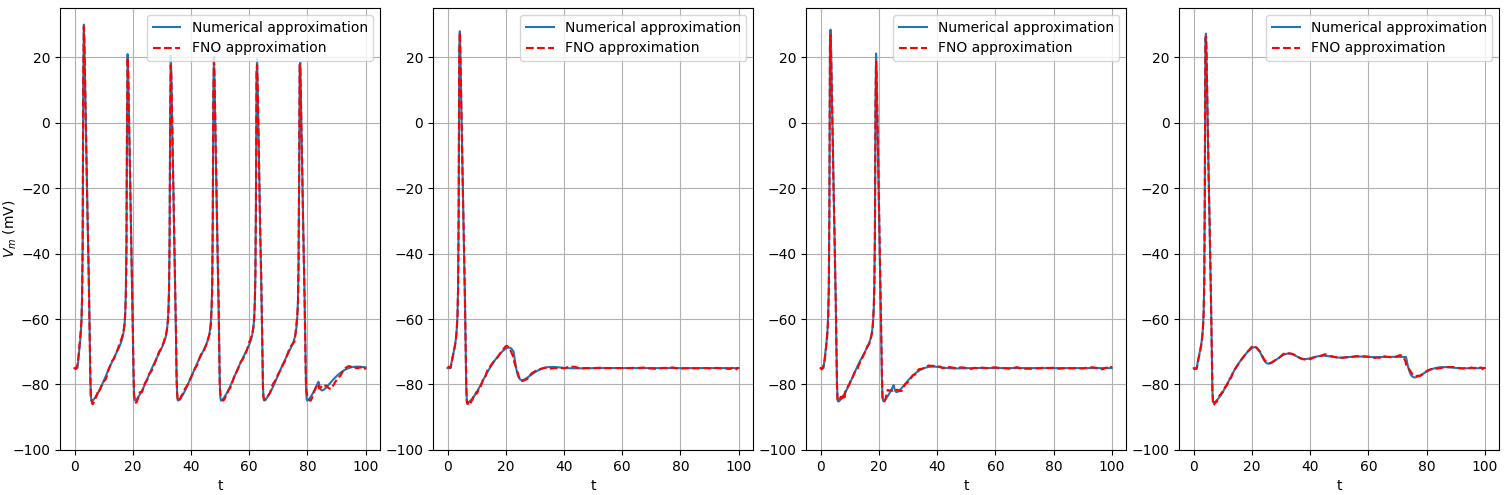}
        \caption{Fourier Neural Operator approximations.}
        \label{fig:FNO_test}
    \end{subfigure}
    \vfill \vspace{0.35cm}
    \begin{subfigure}[b]{\textwidth}
        \centering
        \includegraphics[width=\textwidth]{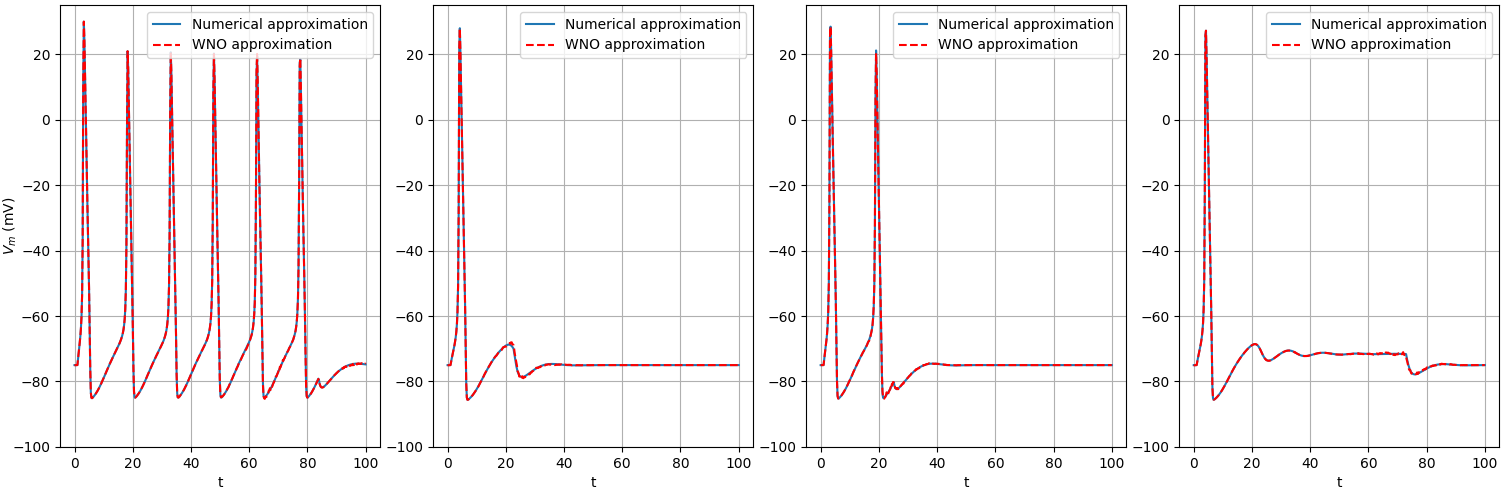}
        \caption{Wavelet Neural Operator approximations.}
        \label{fig:WNO_test}
    \end{subfigure}
    \vfill \vspace{0.35cm}
    \begin{subfigure}[b]{\textwidth}
        \centering
        \includegraphics[width=\textwidth]{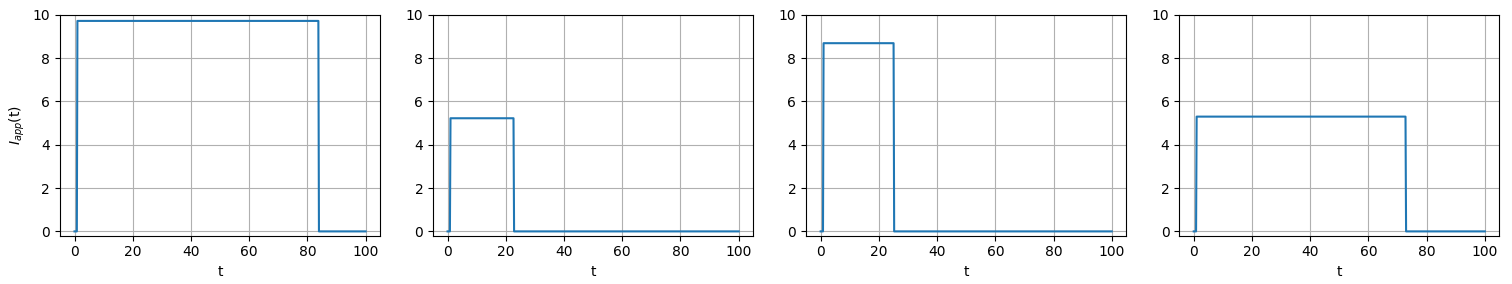}
        \caption{Input functions.}
        \label{fig:input}
    \end{subfigure}
    
    \caption{ (\ref{fig:DON_test}) DeepONets approximations compared to \texttt{ode15s} solutions. (\ref{fig:FNO_test}) FNO approximations compared to \texttt{ode15s} solutions. (\ref{fig:WNO_test}) WNO approximations compared to \texttt{ode15s} solutions. (\ref{fig:input}) Plot of four samples from the test set of the input applied currents functions. In the plots above the numerical solution obtained with \texttt{ode15s} is plotted in blue and the approximation obtained with neural operators is in red dashed line.}
    
    \label{fig:DON_FNO_WNO_input}
\end{figure}

\begin{figure}
    \centering
    \!\!\!\!\includegraphics[width=0.35\textwidth]{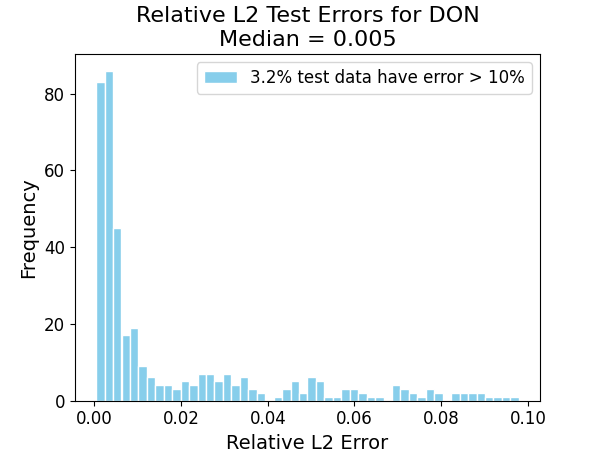}%
    \includegraphics[width=0.35\textwidth]{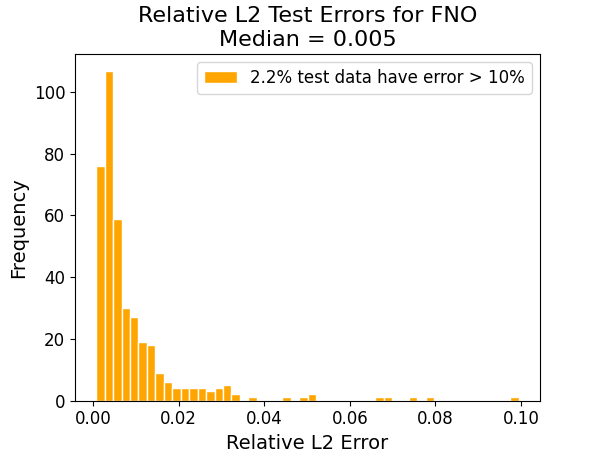}%
    \includegraphics[width=0.35\textwidth]{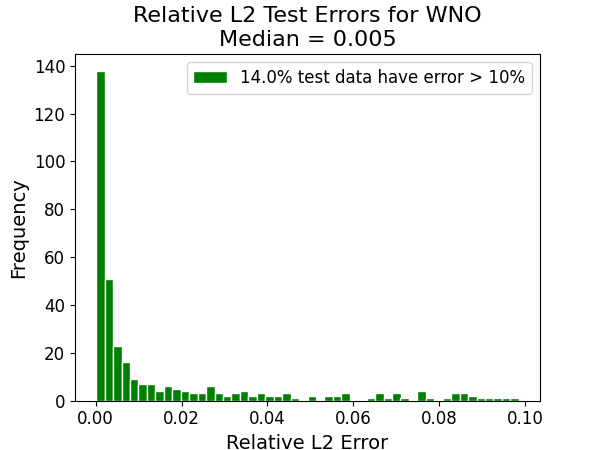}
    \caption{Distribution of the relative $L^2$ validation error on the test dataset each one selected from our best training for each case. FNO has shown less variance around the median validation error value. Despite showing a similar behaviour, WNO trained with the described configuration has 14\% of outliers over 10\% of relative error. DON shows more variance between 0\% and 10\% of relative test error, but not many outliers, which is also the case of FNO. The median value is the same for all the architectures studied.}
    \label{fig:L2histos} 
\end{figure}

\noindent In Fig. \ref{fig:FNO_test} we report four examples of solutions approximated by FNO, for our best model. The four corresponding inputs and the numerical solutions are taken from the test dataset. 
In this case, we observe that both single and many peaks solutions are approximated properly and there is little variance around the mean test error, as also reported in Fig. \ref{fig:L2histos}.

\begin{figure}
    \centering
    \includegraphics[width=0.45\textwidth]{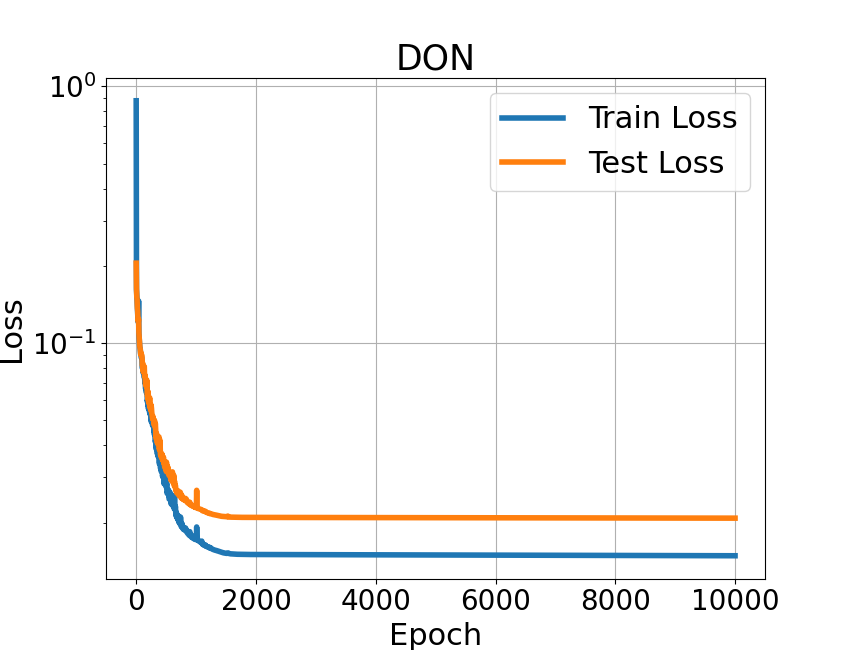}%
    \includegraphics[width=0.45\textwidth]{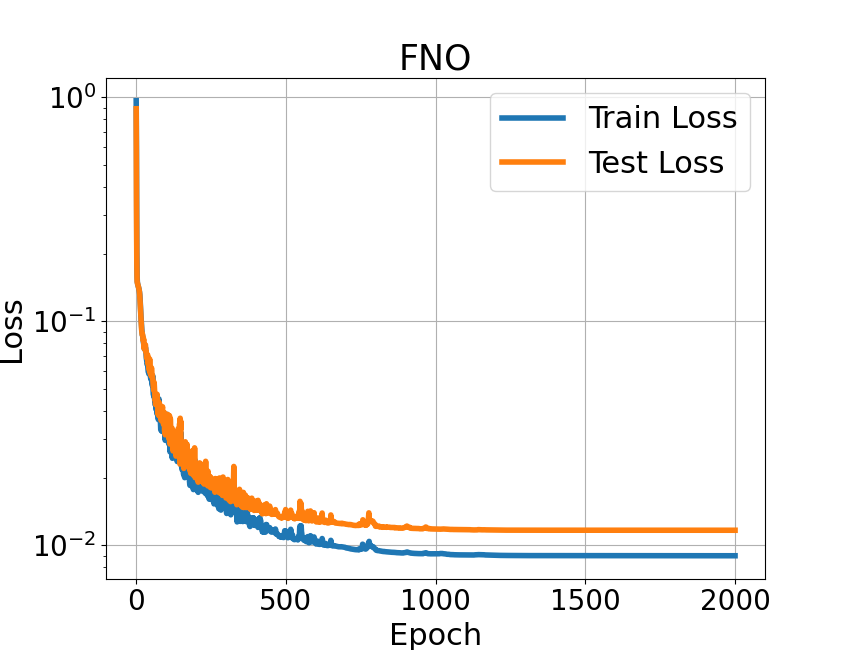}
    \includegraphics[width=0.45\textwidth]{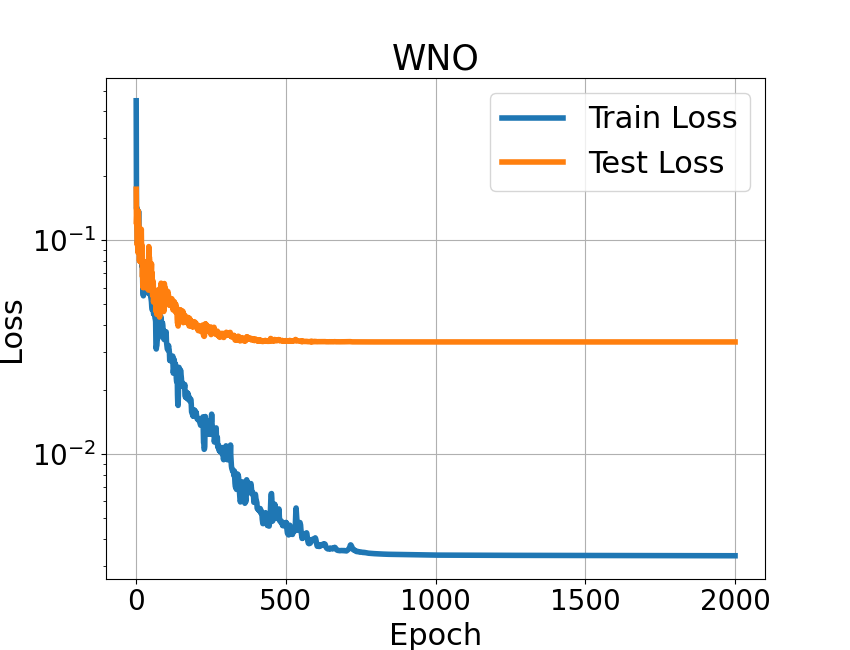}
    \caption{Plot of the train and test relative $L^2$ loss selected from our best training for each case. While DON and FNO show similar train and test errors, WNO trained with the described configuration has shown about an order of magnitude gap between the train and the test error during the training.}
    \label{fig:L2losses} 
\end{figure}

\subsection{WNO results}
\noindent In all the WNO tests performed, we considered all the $500$ discretization points as input to the network, without further encoding.

\noindent In Table \ref{table:results} we have reported our best results for this architecture. Among the hyper-parameters of the best performing trained model we highlight: Gaussian Error Linear Unit (GELU) \cite{hen23} activation function, Daubechies wavelet basis \texttt{db24} for the DWT implemented in the integral kernel, number of integral layers $L=3$ and width of the integral kernel parameterized in the wavelet space $R_{\theta_t}$ equal to $256$. The best mean test loss, evaluated from three trained networks each one with different randomly initialized parameters is around 3.3\% and has been obtained considering a $R_{\theta_t}$ width of 64, for a total of 13 million parameters, the highest number among the architectures considered.

\noindent In Fig. \ref{fig:wno_scatter} we study how the relative $L^2$ error varies when we change the width of the integral kernel parameterized in the wavelet space, $R_{\theta_t}$. We observe a general superlinear reduction of the training error with respect to the width of $R_{\theta_t}$, while the test error appears to saturate after a width of 64. We also observe that, despite having obtained good results in terms of training and test losses (also see Fig. \ref{fig:L2losses}), our best case has worse generalization capabilities and employs a much larger set of parameters, compared to DON and FNO.

\noindent In Fig. \ref{fig:WNO_test} we report four examples of solutions approximated by WNO, for our best model. The four corresponding inputs and the numerical solutions are taken from the test dataset. 
In this case, we observe that both single and many peaks solutions can be approximated properly, but the trained network appears more unstable in keeping the test error around the mean value if compared to DON and FNO, as shown in Fig. \ref{fig:L2histos}.

\subsection{Training Comparison}
\noindent In this section we compare the training of the three different architectures considered.

\noindent In Fig. \ref{fig:L2losses} the evolution of the loss function through the training epochs, namely plots of the relative $L^2$ error (\ref{eq:rel2}) are reported for each one of the architectures object of our study. The FNO has shown less difference between the train and the test error. The DON and the WNO have presented higher loss values in approximating the solutions to our problem, but while the training error in the DON case is similar to the validation one, in the WNO after 50 epochs the training error decays significantly, while the test error has a lower decay rate and saturates quickly. This effect is particularly evident when the width of the $R_{\theta_t}$ operator is increased.

\noindent In Fig. \ref{fig:L2histos} histograms showing the relative $L^2$ error distributions over all test cases for our dataset are reported. An effective training is expected to have an high percentage of the total samples with a relative error around the mean value. For the DON and the FNO all the samples are contained in an interval between $<1\%$ and $16\%$ of validation error. The WNO instead contains all the samples in a wider interval, with a single outlier greater than $25\%$ of relative error.

\section{Conclusions}
\noindent In this work, we have constructed and studied three different operator learning architectures, based on DeepONet, Fourier and Wavelet Neural Operator techniques, to learn the membrane potential dynamics of the Hodgkin-Huxley ionic model. 

\noindent FNOs have shown promising results in predicting the membrane potential when samples are taken from the test dataset, with a mean validation error of 1.4\%. However, it is well known that FNOs need an equispaced domain in order to work properly, which means a potentially large number of domain points when facing stiff problems. DeepONets are more flexible in this sense, allowing adaptive and irregular domains. However, for the Hodgkin-Huxley model they have shown less capability to learn many-peaks solutions, compared to the other architectures. Their mean validation error of 2.2\% with our best configuration is still an important result. For further studies it can be useful to explore for this class of problems also different DeepONet based architectures, e.g. S-DeepONets \cite{hex24}. WNOs have shown a very good capability of generalizing solutions with many peaks over the test dataset, but they required more tuning in the number of integral layers and in the wavelet type employed, as well as a larger number of trainable parameters and a consequent larger training time. Also, the mean validation error of 3.3\% is the highest among the architectures studied.

\noindent As a final remark, while at this level there is no real need of employing these methods in order to get a more efficient solver in terms of time or complexity for the Hodgkin-Huxley model if compared to the classical ones (e.g. Runge-Kutta), it is important to study systems which can be challenging to learn by artificial neural networks due to their properties, in particular the stiffness of the solutions and the threshold-dependent dynamics. This work paves the way for employing the same techniques for a whole family of similar systems, with way more equations involved or even with more dimensions: further studies in this directions can be useful to design architectures capable of solving more physiological models in computational cardiology or neuroscience, making properly trained operator learning models a valid aid for personalized medicine. 

\section{Acknowledgements}
\noindent The authors wish to thank Carlo Marcati for helpful discussions and suggestions. Edoardo Centofanti and Luca F. Pavarino have been supported by grants of INdAM–GNCS and MUR (PRIN 202232A8AN\_002 and PRIN P2022B38NR\_001),  funded by European Union - Next Generation EU.

\section{Code Availability}\label{sec:avail}
\noindent
The codes for the numerical tests presented in this work will be made publicly available on GitHub.

\end{document}